# Self-morphing of elastic bilayers induced by mismatch strain: deformation simulation and bio-inspired design

*Junjie Song[1,2], Yixiong Feng[1,2]\*, Zhaoxi Hong[1,2], Bingtao Hu[1,2], and Jianrong Tan[1,2], Xiuju Song[1,2]\**

1 State Key Laboratory of Fluid Power and Mechatronic Systems, Zhejiang University, Hangzhou 310027, P. R. China

2 Engineering Research Center for Design Engineering and Digital Twin of Zhejiang Province, Zhejiang University, Hangzhou 310027, P. R. China

E-mail: songxiuju@zju.edu.cn

**Abstract**

The process of self-morphing in curved surfaces found in nature, such as with the growth of flowers and leaves, has generated interest in the study of self-morphing bilayers, which has been used in many soft robots or switchers. However, previous research has primarily focused on materials or bilayer fabrication technologies. The self-morphing mechanism and process have been rarely investigated, despite their importance. This study proposed a new deformation simulation method for self-morphing bilayers based on a checkerboard-based discrete differential geometry approach. This new method achieved higher efficiency than traditional finite element methods while still maintaining accuracy. It was also effective in handling complex finite strain situations. Finally, the simulation model was used to design three self-morphing bilayers inspired by folding flowers, spiral grass, and conical seashells. These designs further prove the effectiveness of the proposed method. The results of this study propose a good method for predicting deformation and designing self-morphing bilayers and provide a useful viewpoint for using geometrical methods to solve mechanical problems.

# 1 Introduction

There is a natural self-morphing process in the growth of flowers, leaves, and other plants where curved surfaces develop. Researchers have been inspired by these remarkable natural surfaces to create artificial self-morphing bilayers that possess similar properties [1-8]. In recent years, a vast range of self-morphing bilayers has been proposed for numerous applications including but not limited to soft grippers [9-13], switchers [14-16], and crawling robots [17-23]. With the development of additive manufacturing (AM) technology and smart materials, people find that combining these two technologies will bring great convenience for the fabrication of these self-morphing bilayers. This led to the emergence of 4D printing technology, where the "4D" refers to the self-morphing behavior of printed structures over time [24]. Then numerous studies were conducted to investigate viable printing technologies (such as fused deposition modeling (FDM) [25-27], digital light processing (DLP) [28,29], and direct ink writing (DIW) [30,31]) and intelligent materials (such as shape memory polymers (SMP) [32,33], hydrogels [34,35], etc.) for 4D printing. However, only a few studies have focused on the self-morphing bilayers' deformation mechanisms or morphology control.

Mismatch strain is the primary driver of self-morphing processes in both natural and artificial bilayers. Specifically, the growth of different layers in these structures occurs in various directions or with varying expansive ratios, which place restrictions on and influence one another, ultimately resulting in coupled deformation. Understanding how mismatch strain impacts bilayer coupling deformation is crucial for controlling self-morphing bilayers' deformation to meet application demands. For simple situations such as bending, deformation can be easily modeled using curvature, and the deformation mechanism has been extensively studied [36-39]. However, explicit relationships for other complex scenarios are challenging to determine. The most common approach is to deduce the bilayers' elastic potential energy function and minimize it as the self-morphing process follows the principle of minimum potential energy [40-44]. Minimizing the energy function can help identify the bilayers' equilibrium morphology, thereby enabling the deformation simulation of self-morphing bilayers.

The elastic potential energy function of the bilayers that is deduced often takes the form of

$$E = \int_{\Omega_h} W(\nabla v) \tag{1}$$

where $v: \Omega_h = S \times \left(-\frac{h}{2}, \frac{h}{2}\right) \to R^3$, $S$ denotes the two-dimensional (2D) domain of the bilayer's mid-surface and $h$ is the thickness of the bilayer. The energy function of the bilayer is a function of the deformation gradient of the mid-surface of the bilayer, which includes the deformation gradient $\nabla v$ and its gradients up to a certain order. The subsequent question is how to solve the minimization problem. Directly minimizing the energy function is challenging since it involves derivatives with different orders, leading to a singular perturbation problem [45]. One approach is to scale Eq. (1) to be $\lim_{h \to 0} \frac{1}{h^3} E$ and minimizing this scaled energy function [40-43]. This operation is capturing the energy term in Eq. (1) that is proportional to $h^3$, corresponding to a bending deformation of the mid-surface [46]. Then traditional methods like calculus of variations and finite element method (FEM) can be introduced to minimize the scaled energy function. However, the drawback of this approach is that it assumes the mid-surface is unstretched during deformation to keep the energy function bounded, which contradicts reality and thus is only suitable for small strain situations (please refer to section 2 for a more detailed explanation). Another way to deal with the energy minimization problem is to discretize Eq. (1) based on discrete differential geometry (DDG) approaches. For instance, Ree et al [44] simplified Eq. (1) as the function of the first and second fundamental quantities of the bilayer's mid-surface and used a DDG approach suitable for general triangle meshes [47] to discretize the energy function. Later, the energy minimization problem was converted to a discrete least square optimization problem solvable through gradient descent algorithms. The main difficulties of this means are how to simplify Eq. (1) into a form that enables easy discretization and how to select a suitable DDG approach with good convergence and high accuracy for discretization.

In 2022, a new approach called the checkerboard-based DDG approach was introduced by Dellinger [48]. The approach's origin can be traced back to Peng's work in 2019 [49], where they initially discovered that isometric (or conformal) mapping between two surfaces could be constructed by congruence (or similarity) constraints between the black inscribed rectangles in the surfaces' checkerboard discretization. Later, Jiang et al. [50,51] extended this idea to the general case by validating the conclusion still holds when removing the rectangle restriction to the black inscribed quadrilaterals. Motivated by the promising outcomes of this approach in surface mapping, Dellinger established the DDG approach based on checkerboard patterns and proved its convergence in theory. In comparison to other curvature discrete methods, this approach has simpler discrete forms and has proven convergence. However, Dellinger's approach is only suitable for structured quadrilateral meshes, whereas unstructured quadrilateral meshes are prevalent for complex surface meshing. To date, the approach has only been used for discrete mapping between surfaces. Nevertheless, it is also suitable for the energy minimization problems of deforming surfaces.

In this study, we derived a simpler elastic energy function for the bilayer structure that is easier to discretize and can reflect the bilayer deformation trend during self-morphing more intuitively. The checkerboard-based DDG approach was improved for unstructured quadrilateral meshes by introducing the mean weighted equaled (MWE) algorithm to calculate vertex normal [52]. Next, we discretized the elastic energy function with the checkerboard-based DDG approach, converting the elastic energy minimization problem of self-morphing bilayers into a least square optimization problem that we solved with the Levenberg-Marquardt (LM) algorithm [53]. This enabled us to establish the checkerboard-based deformation simulation model of self-morphing bilayers. Then, we conducted several numerical experiments on the bending deformation of rectangle self-morphing bilayers, including both small strain situations and finite strain situations. The simulation results of the newly proposed simulation model were compared with the results of a FEM-based simulation model (proposed in our previous work [42]), to verify the high efficiency and accuracy of the new model. Additionally, the numerical experiments on finite strain situations demonstrated the superiority of the new model, which considers the bilayer's stretch during deformation. Finally, we designed three kinds of self-morphing bilayers inspired by four-petaled folding flowers, spiral grasses, and conical seashells, based on the proposed model, and fabricated them using FDM-based 4D printing technology. The congruity between simulation results and experiments further certified the proposed model's effectiveness.

This paper is organized as follows: In Section 2, we deduced the simplified form of the bilayers' elastic energy function. In Section 3, we improved the checkerboard-based DDG approach for the general case and utilized it to discretize the energy function derived in Section 2. Subsequently, we established the deformation simulation model for self-morphing bilayers. In Section 4, we conducted numerical experiments on bending rectangle bilayers using both the proposed model and the FEM-based model in our previous work; we compared the simulation outcomes of the two models. In Section 5, we designed three bio-inspired self-morphing bilayers with 4D printing technology based on the proposed model and compared simulation outcomes with experiments. Finally, in Section 6, we drew some conclusions and indicated possible future directions.

**2 Deformation modeling of a self-morphing bilayer**

As mentioned before, the deformation of a self-morphing bilayer obeys the principle of minimum potential energy. The deformation modeling and simulation of a self-morphing bilayer equals to deduce the elastic energy function of a self-morphing bilayer and realize the energy minimization.

To obtain the elastic energy function of a self-morphing bilayer, two steps are necessary: deducing the bilayer's strain tensor and simplifying the strain energy function. Given the bilayer's small thickness, we introduce the assumptions of Kirchhoff-Love plate theory in our study to simplify the bilayer's deformation modeling.

*Kirchhoff-Love assumptions:* a) straight lines normal to the mid-surface remain straight after deformation; b) straight lines normal to the mid-surface remain normal to the mid-surface after deformation; c) the thickness of the plate does not change during deformation.

The self-morphing bilayer consists of two distinct layers with different strains. To enhance clarity regarding the bilayer plate's deformation process, we define three states (as illustrated in Fig. 1) for a monolayer within the bilayer:

1. The Original State: This state signifies the initial shape of the monolayer before any deformation. It serves as a starting point for the deformation analysis.

2. The Relaxed State: This state corresponds to an imaginary state of the monolayer and denotes the deformed shape of the monolayer when it is not subject to any constraints from the adjacent monolayer. It serves as a reference point to comprehend the deformation caused solely by internal growth factors, excluding the influence of the neighboring layer.

3. The Current State: This state refers to the deformation of the monolayer caused by the strain mismatch between the two layers in the bilayer. It represents the bilayer's shape after undergoing deformation.

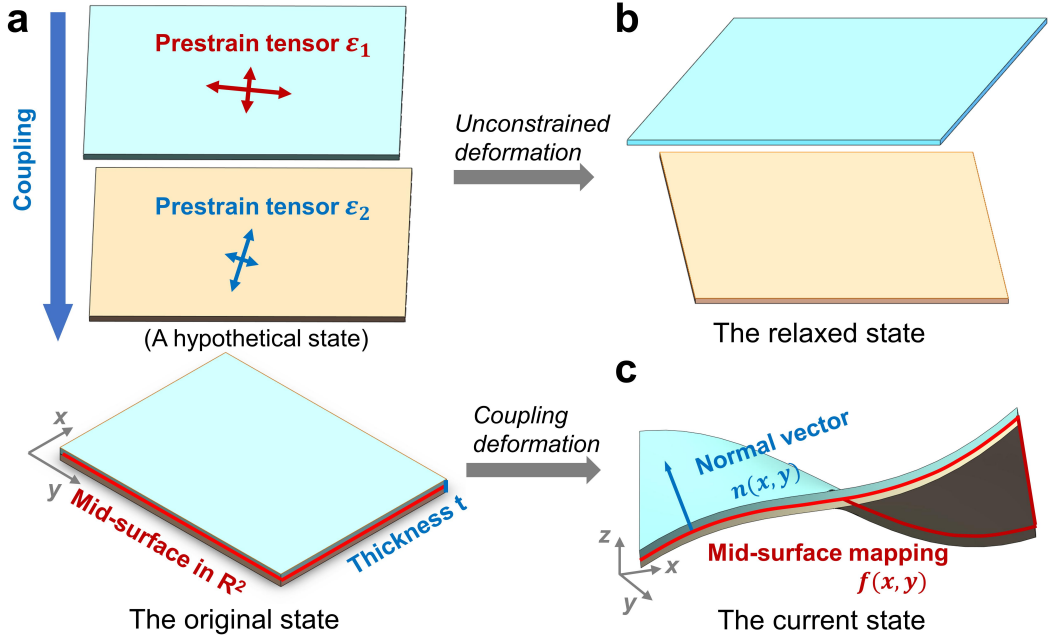

**Fig. 1.** Three states of a self-morphing bilayer during the deformation process: (a) The original state: the initial shape of the monolayer before any deformation. (b) The relaxed state: the unconstrained deformation of monolayers in the bilayer (an imaginary state). (c) The current state: the coupling shape after undergoing deformation.

2.1 Finite strain tensor of the self-morphing bilayer

The deformation of the monolayer within the bilayer can be considered a two-step process. The first step involves the deformation of the monolayer from its relaxed state to its original state. The second step involves the deformation from the original state to the current state with prestrain resulting from the first step. As a result, the finite strain tensor, which takes into account the prestrain of the bilayer, can be expressed as follows:

$$\boldsymbol{\epsilon} = \frac{1}{2}(F_c^T F_c - F_r^T F_r) \qquad (2)$$

where $F_r$ and $F_c$ are the deformation gradient tensor of the relaxed state and the current state respectively, with the original state as reference.

Here we define a linear strain tensor, following to the Kirchhoff-Love assumptions, to

describe the unconstrained in-plane growth of a monolayer, which can be expressed as

$$\boldsymbol{\varepsilon} = \begin{bmatrix} \varepsilon_{11} & \varepsilon_{12} \\ \varepsilon_{21} & \varepsilon_{22} \end{bmatrix} \tag{3}$$

Then the deformation gradient tensor $F_r$ can be calculated as

$$F_r = I_3 + \begin{bmatrix} \boldsymbol{\varepsilon} & 0 \\ 0 & 1 \end{bmatrix} = \begin{bmatrix} 1+\varepsilon_{11} & \varepsilon_{12} & 0 \\ \varepsilon_{21} & 1+\varepsilon_{22} & 0 \\ 0 & 0 & 1 \end{bmatrix} \tag{4}$$

The deduction of the deformation gradient $F_c$ needs to be conducted based on the parametric model of the bilayer. The deformation of the bilayer, denoted as $g: R^3 \to R^3$, can be approximated under the Kirchhoff-Love assumptions

$$g(x, y, z) = f(x, y) + z \cdot n(x, y) \tag{5}$$

where $(x, y)$ are the parameters defining the original state, $f: R^2 \to R^3$ is the deformation of the bilayer mid-surface, $n(x, y)$ is the normal vectors of the surface mapping $f(x, y)$, $z \in [-t/2, t/2]$ represents the z coordinate along the thickness of the bilayer, as shown in Fig. 1. The deformation gradient $F_c$ can be calculated based on the mapping $g(x, y, z)$ as

$$F_c = \begin{bmatrix} \frac{\partial f}{\partial x} + z\frac{\partial n}{\partial x} & \frac{\partial f}{\partial y} + z\frac{\partial n}{\partial y} & n \end{bmatrix} \tag{6}$$

2.2 Elastic energy function of the self-morphing bilayer

In this study, the simplified strain energy density function in [41]

$$W(\epsilon) = \mu|\epsilon|^2 \tag{7}$$

is introduced to describe the elastic energy of the bilayer, where $\mu$ is the second *Lamé* constant and $|H|$ is the Frobenius norm of the matrix $H$. The elastic energy function of the self-morphing bilayer is obtained by integrating the finite strain tensor over the entire bilayer. Assuming that the two layers of the bilayer have different second Lamé constants, represented by $\mu_1$ and $\mu_2$, different thicknesses, represented by $t_1$ and $t_2$, and different linear strain tensors, represented by $\varepsilon_1$ and $\varepsilon_2$, the elastic energy function can be expressed by combining Eqs. (2), (4), (6), and (7),

as follows:

$$E[f] = \int W(\epsilon) = \int_\omega \int_t \frac{\mu}{4} \left\| \begin{bmatrix} (df+zdn)^T(df+zdn) & 0 \\ 0 & n^Tn \end{bmatrix} - \begin{bmatrix} (I_2+\epsilon)^2 & 0 \\ 0 & 1 \end{bmatrix} \right\|^2$$

$$= \int_\omega \int_t \frac{\mu}{4} |a_g - 2zb_g + z^2 c_g - I_2 - 2\epsilon - \epsilon^2|^2 \qquad (8)$$

where $a_g$ and $b_g$ are the first and second fundamental forms of the surface mapping $f$. We can neglect the term $z^2 c_g$ and expand the expression in the integral, as $c_g = dn^T dn$ and the thickness of the bilayer is small compared to other terms

$$E[f] = \int_\omega \int_t \frac{\mu}{4} |a_g - 2zb_g - I_2 - 2\epsilon - \epsilon^2|^2$$

$$= \int_\omega \int_t \frac{\mu}{4} | a_g^2 + 4b_g^2 z^2 + I_2 + 4\epsilon^2 + \epsilon^4 - 4a_g b_g z - 2a_g - 4\epsilon a_g - 2\epsilon^2 a_g + 4b_g z + 8\epsilon b_g z$$

$$+ 4\epsilon^2 b_g z + 4\epsilon + 2\epsilon^2 + 4\epsilon^3 |$$

$$= \int_\omega \int_{\frac{t_2-t_1}{2}}^{\frac{t_1+t_2}{2}} \frac{\mu_1}{4} |4b_g^2 z^2 + (4\epsilon_1^2 b_g + 8\epsilon_1 b_g + 4b_g - 4a_g b_g)z + a_g^2 + I_2 + 4\epsilon_1^2 + \epsilon_1^4 - 2a_g -$$

$$4\epsilon_1 a_g - 2\epsilon_1^2 a_g + 4\epsilon_1 + 2\epsilon_1^2 + 4\epsilon_1^3 | + \int_\omega \int_{-\frac{t_1^2+t_2}{2}}^{\frac{t_2-t_1}{2}} \frac{\mu_2}{4} |4b_g^2 z^2 + (4\epsilon_2^2 b_g + 8\epsilon_2 b_g + 4b_g -$$

$$4a_g b_g)z + a_g^2 + I_2 + 4\epsilon_2^2 + \epsilon_2^4 - 2a_g - 4\epsilon_2 a_g - 2\epsilon_2^2 a_g + 4\epsilon_2 + 2\epsilon_2^2 + 4\epsilon_2^3 |$$

$$= \int_\omega \frac{(\mu_1 t_1 + \mu_2 t_2)}{4} a_g^2 - \frac{\mu_1 t_1 (I_2+\epsilon_1)^2 + \mu_2 t_2 (I_2+\epsilon_2)^2}{2} a_g + \frac{\mu_1 t_1 (3t_2^2+t_1^2) + \mu_2 t_2 (3t_1^2+t_2^2)}{3} b_g^2 + 2t_1 t_2 (\mu_1 \epsilon_1^2 +$$

$$2\mu_1 \epsilon_1 - \mu_2 \epsilon_2^2 - 2\mu_2 \epsilon_2) b_g + C_1$$

$$= \int_\omega \frac{(\mu_1 t_1 + \mu_2 t_2)}{4} \left| a_g - \frac{\mu_1 t_1 (I_2+\epsilon_1)^2 + \mu_2 t_2 (I_2+\epsilon_2)^2}{(\mu_1 t_1 + \mu_2 t_2)} \right|^2 + \frac{\mu_1 t_1 (3t_2^2+t_1^2) + \mu_2 t_2 (3t_1^2+t_2^2)}{12} \Big| b_g -$$

$$\frac{3t_1 t_2 (\mu_2 \epsilon_2^2 + 2\mu_2 \epsilon_2 - \mu_1 \epsilon_1^2 - 2\mu_1 \epsilon_1)}{\mu_1 t_1 (3t_2^2+t_1^2) + \mu_2 t_2 (3t_1^2+t_2^2)} \Big|^2 + C_2 \quad (9)$$

Here $C_1$, $C_2$ are constants that have no effects on the minimizing process of the energy function.

Finally, the elastic energy function of the self-morphing bilayer can take the form of

$$\widetilde{E}[f] := \int_\omega \omega_{first} |a_g - A|^2 + \omega_{second} |b_g - B|^2 \qquad (10)$$

$$\text{where } \omega_{first} = \frac{(\mu_1 t_1 + \mu_2 t_2)}{4}, \qquad (11)$$

$$\omega_{second} = \frac{\mu_1 t_1 (3t_2^2+t_1^2) + \mu_2 t_2 (3t_1^2+t_2^2)}{12}, \qquad (12)$$

$$A = \frac{\mu_1 t_1 (I_2+\epsilon_1)^2 + \mu_2 t_2 (I_2+\epsilon_2)^2}{(\mu_1 t_1 + \mu_2 t_2)}, \qquad (13)$$

$$B = \frac{3t_1 t_2 (\mu_2 \varepsilon_2^2 + 2\mu_2 \varepsilon_2 - \mu_1 \varepsilon_1^2 - 2\mu_1 \varepsilon_1)}{\mu_1 t_1 (3t_2^2 + t_1^2) + \mu_2 t_2 (3t_1^2 + t_2^2)}. \tag{14}$$

Equation (10) intuitively showed that during the self-morphing, the stretching and bending deformations of the bilayer both occur, corresponding to the changes of the first and second fundamental quantities of the bilayer, respectively. In approaches of scaling the energy function [40-43], the term $A$ was considered to be $I_2$ and $a_g - A$ was assumed to be zero to keep $\widetilde{E}[f]/h^3$ bounded under the scaled factor $1/h^3$. This allowed the minimization of the rest energy term $\frac{1}{h^3}\int_\omega \omega_{second}|b_g - B|^2$ to be solvable. However, the isometric constraint, which assumption $a_g = I_2$, only holds for small strain situations. This assumption leads to significant errors when the strain is not small (see Section 4). One of the advantages of the proposed model in this study is that it considers both stretching and bending deformations of self-morphing bilayers during simulation, ensuring good accuracy even in situations with significant strain.

**3 Deformation simulation of a self-morphing bilayer based on checkerboard patterns**

Here we introduced a new proposed DDG approach, the checkerboard patterns [48], to discretize the deduced energy function Eq. (10). Hence the minimization of Eq. (10) can be converted to a nonlinear least square optimization problem, which greatly reduced the difficulty of problem solving.

To obtain the checkerboard patterns of a quadrilateral mesh, the midpoints of each edge of every quadrilateral are connected in order and the inscribed parallelogram is obtained, as illustrated in Fig. 2a.

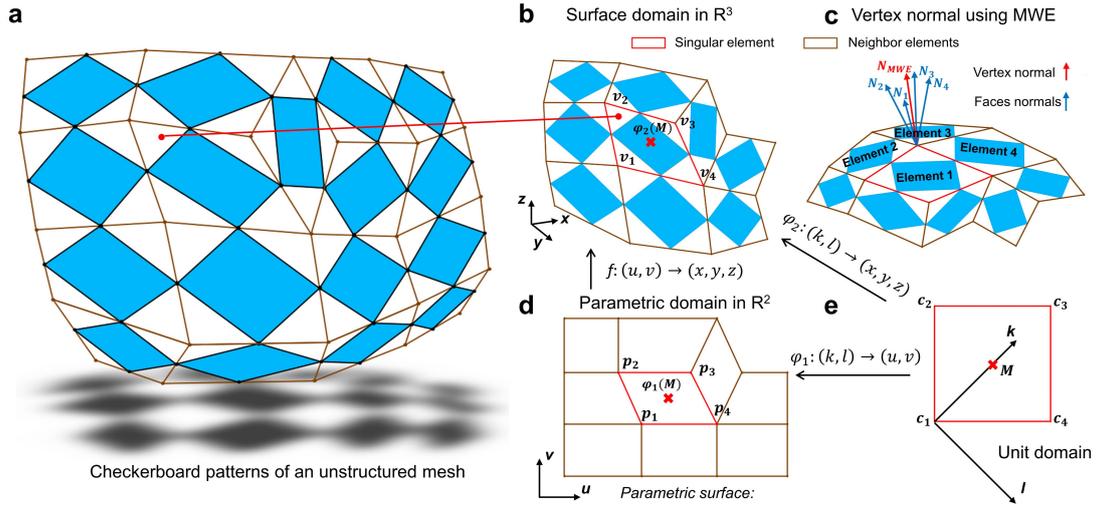

**Fig. 2.** (a) The checkerboard patterns of an unstructured 3D surface mesh. (b) A singular element and its neighbor elements on the surface in $R^3$. (c) The illustration of the MWE algorithm used to calculate the vertex normal. (d) The mapping elements of the singular element and its neighbor elements in parametric domain in $R^2$. (e) The $1 \times 1$ unit domain with the basis $k = \frac{1}{\sqrt{2}}(1,1)^T$ and $l = \frac{1}{\sqrt{2}}(1,-1)^T$ in $R^2$.

Considering a parametric surface $f: R^2 \to R^3$, mesh the parametric domain and surface domain with a pair of combinatorially equivalent unstructured quadrilateral meshes. To illustrate the discrete process based on checkerboard patterns, we study an arbitrary singular element in the unstructured meshes (where "singular" means the degrees of the vertexes in this element are not all equal to 4).

3.1 Discretization of the first fundamental form

Figure 2d shows a singular element and its neighboring elements in $R^2$ (in the parametric domain), while Fig. 2b shows the corresponding mapping elements in $R^3$ (on the surface). To facilitate analysis, the vertices of the singular element in the meshes are numbered sequentially, as seen in Figure 2. We also construct a $1 \times 1$ unit domain $c_1 c_2 c_3 c_4$, as shown in Fig. 2e, and establish the mappings from this domain to the singular elements in $R^2$ and $R^3$ as $\varphi_1: R^2 \to R^2$ and $\varphi_2: R^2 \to R^3$, respectively (see Fig. 2).

Using the basis $k = \frac{1}{\sqrt{2}}(1,1)^T$ and $l = \frac{1}{\sqrt{2}}(1,-1)^T$ to define a coordinate system of the unit square, the derivatives $\partial_k \varphi_2$ and $\partial_l \varphi_2$ at the center point M (see Fig. 2e) can be approximated by central difference as

$$\partial_k \varphi_2(M) \approx \frac{\varphi_2(c_3) - \varphi_2(c_1)}{\sqrt{2}} = \frac{v_3 - v_1}{\sqrt{2}} \tag{15}$$

$$\partial_l \varphi_2(M) \approx \frac{\varphi_2(c_4) - \varphi_2(c_2)}{\sqrt{2}} = \frac{v_4 - v_2}{\sqrt{2}} \tag{16}$$

where $v_i$ refers the coordinates of the $i$th vertex of the singular element in $R^3$. Then the first fundamental form of the parametric surface under the parameters $(k, l)$ can be calculated as

$$I_{(k,l)} = \begin{bmatrix} \partial_k \varphi_2 \cdot \partial_k \varphi_2 & \partial_k \varphi_2 \cdot \partial_l \varphi_2 \\ \partial_l \varphi_2 \cdot \partial_k \varphi_2 & \partial_l \varphi_2 \cdot \partial_l \varphi_2 \end{bmatrix} \tag{17}$$

For parametric surfaces, the fundamental forms under the parameters $(u, v)$ are always of more concern, where, in general, $u = (1,0)^T$ and $v = (0,1)^T$ define the coordinate system of the parametric domain (see Fig. 2d). To obtain the fundamental forms under $(u, v)$, the Jacobian matrix of the mapping $\varphi_1: R^2 \to R^2$ is needed. By introducing the three-dimensional (3D) homogeneous space, the 2D mapping $\varphi_1$ can be expressed as

$$\begin{bmatrix} a & d & g \\ b & e & h \\ c & f & 1 \end{bmatrix} \begin{bmatrix} k \\ l \\ 1 \end{bmatrix} = \begin{bmatrix} wu \\ wv \\ w \end{bmatrix} \to \begin{cases} u = \frac{ak+bl+c}{gk+hl+1} \\ v = \frac{dk+el+f}{gk+hl+1} \end{cases} \tag{18}$$

where $[k \ \ l \ \ 1]^T$ is the homogeneous coordinates of the unit domain, $[wu \ \ wv \ \ w]^T$ is the homogeneous coordinates of the parametric domain ($w$ can be any non-zero value) and the $3 \times 3$ matrix is the homogeneous transform matrix representing the mapping $\varphi_1$.

By substituting the coordinates of points $c_1(0,0)$, $c_2(\sqrt{2}/2, -\sqrt{2}/2)$, $c_3(\sqrt{2}, 0)$, $c_4(\sqrt{2}/2, \sqrt{2}/2)$ and $p_1(u_1, v_1)$, $p_2(u_2, v_2)$, $p_3(u_3, v_3)$, $p_4(u_4, v_4)$ (see Fig. 2), we can obtain the linear equations to solve the parameters $a \sim g$ in the homogeneous transform matrix as

$$\begin{bmatrix} 1 & -1 & 0 & 0 & 0 & 0 & -u_2 & u_2 \\ 0 & 0 & 0 & 1 & -1 & 0 & -v_2 & v_2 \\ 0 & 0 & 1 & 0 & 0 & 0 & 0 & 0 \\ 1 & 0 & 0 & 0 & 0 & 0 & -u_3 & 0 \\ 0 & 0 & 0 & 1 & 0 & 0 & -v_3 & 0 \\ 0 & 0 & 0 & 0 & 0 & 1 & 0 & 0 \\ 1 & 1 & 0 & 0 & 0 & 0 & -u_4 & -u_4 \\ 0 & 0 & 0 & 1 & 1 & 0 & -v_4 & -v_4 \end{bmatrix} \begin{bmatrix} a \\ b \\ c \\ d \\ e \\ f \\ g \\ h \end{bmatrix} = \begin{bmatrix} \sqrt{2}(u_2 - u_1) \\ \sqrt{2}(v_2 - v_1) \\ u_1 \\ (u_3 - u_1)/\sqrt{2} \\ (v_3 - v_1)/\sqrt{2} \\ v_1 \\ \sqrt{2}(u_4 - u_1) \\ \sqrt{2}(v_4 - v_1) \end{bmatrix} \quad (19)$$

Then the Jacobian matrix of the mapping $\varphi_1$ can be calculated with the solved parameters as

$$J = \begin{bmatrix} \frac{\partial u}{\partial k} & \frac{\partial u}{\partial l} \\ \frac{\partial v}{\partial k} & \frac{\partial v}{\partial l} \end{bmatrix} = \begin{bmatrix} \frac{(ah-gb)l+a-gc}{(gk+hl+1)^2} & \frac{(gb-ah)k+b-hc}{(gk+hl+1)^2} \\ \frac{(dh-ge)l+d-gf}{(gk+hl+1)^2} & \frac{(eg-dh)k+e-hf}{(gk+hl+1)^2} \end{bmatrix} \quad (20)$$

As $f: \varphi_1^{-1} \circ \varphi_2$, the first fundamental form of the parametric surface under $(u, v)$, at the center point M, is deduced as

$$I_{(u,v)} = (J^{-1})^T \begin{bmatrix} \partial_k \varphi_2 \cdot \partial_k \varphi_2 & \partial_k \varphi_2 \cdot \partial_l \varphi_2 \\ \partial_l \varphi_2 \cdot \partial_k \varphi_2 & \partial_l \varphi_2 \cdot \partial_l \varphi_2 \end{bmatrix} J^{-1} \quad (21)$$

Then the first fundamental form of any element in the parametric surface can be obtained with Eq. (21).

3.2 Discretization of the second fundamental form

To discretize the second fundamental form, we need to first build the Gauss map of the parametric surface, which is equivalent to calculating the surface normal. In Dellinger's theory, the surface normal was defined at every mesh vertex and calculated using the cross-product of the coordinates of the four neighboring vertices. However, this method is limited to quadrilateral meshes where all vertices have four degrees.

To extend this method to unstructured quadrilateral meshes, we define the surface normal at every mesh vertex using the MWE algorithm. This algorithm efficiently computes vertex normals and has discrepancies of no more than 1-2° in most cases. Moreover, these discrepancies decrease as the mesh is refined (see [52]).

*The MWE algorithm* the normal vector on a vertex is considered to be parallel to the

summation of all face normals of the faces incident to the vertex:

$$N_{MWE} \parallel \sum_{i=1}^{n} N_i \tag{22}$$

where $n$ is the number of faces incident to the vertex and $N_i$ is the face normal of the $i$th face, see Fig. 2c.

We use the cross product of the two diagonal vectors of a quadrilateral to calculate the face normal, taking the face normal of the singular quadrilateral as an example:

$$N_s = \frac{(v_3-v_1)\times(v_4-v_2)}{|(v_3-v_1)\times(v_4-v_2)|} \tag{23}$$

Then we can use the same way like Eqs. (15) and (16) to approximate the derivatives of surface normal on the singular quadrilateral as

$$\partial_k n \approx \frac{n_3-n_1}{\sqrt{2}} \tag{24}$$

$$\partial_l n \approx \frac{n_4-n_2}{\sqrt{2}} \tag{25}$$

where $n$ is the surface normal at the center point of the singular quadrilateral, $n_i$ is the vertex normal of the $i$th vertex. Thus, the second fundamental form of the parametric surface under the parameters $(k, l)$ can be expressed as

$$II_{(k,l)} = \begin{bmatrix} \partial_k n \cdot \partial_k \varphi_2 & \partial_k n \cdot \partial_l \varphi_2 \\ \partial_l n \cdot \partial_k \varphi_2 & \partial_l n \cdot \partial_l \varphi_2 \end{bmatrix} \tag{26}$$

Similarly, we can calculate the second fundamental form of the parametric surface under $(u, v)$ as:

$$II_{(u,v)} = (J^{-1})^T \begin{bmatrix} \partial_k n \cdot \partial_k \varphi_2 & \partial_k n \cdot \partial_l \varphi_2 \\ \partial_l n \cdot \partial_k \varphi_2 & \partial_l n \cdot \partial_l \varphi_2 \end{bmatrix} J^{-1} \tag{27}$$

3.3 Energy minimization via non-linear least squares optimization

A 2D quadrilateral mesh $M = (V, E, F)$ is used to discretize the bilayer mid-surface domain as the initial mesh before deformation, where V, E, and F are the vertex, edge and face sets of the mesh. Based on the checkerboard-based DDG approach, the minimization of Eq. (10) can be converted to a non-linear least squares problem with the form as

$$E_{discrete} = \omega_{first} \sum_{f \in F} \sum_{j=1}^{3} c_{first,j}(f)^2 + \omega_{second} \sum_{f \in F} \sum_{j=1}^{4} c_{second,j}(f)^2 \quad (28)$$

The optimization objective functions correspond to the first and second forms of the surface. Considering a face in the mesh $f = v_1 v_2 v_3 v_4$ (the numbering sequence is the same as Fig. 2b), the optimization objective functions can be expressed as follows:

$$c_{first,1}(f) = \frac{(v_3 - v_1) \cdot (v_3 - v_1)}{2} - \bar{A}_{11} \quad (29)$$

$$c_{first,2}(f) = \frac{\sqrt{2}(v_3 - v_1) \cdot (v_4 - v_2)}{2} - \sqrt{2}\bar{A}_{12} \quad (30)$$

$$c_{first,3}(f) = \frac{(v_4 - v_2) \cdot (v_4 - v_2)}{2} - \bar{A}_{22} \quad (31)$$

$$c_{second,1}(f) = \frac{(n_3 - n_1) \cdot (v_3 - v_1)}{2} - \bar{B}_{11} \quad (32)$$

$$c_{second,2}(f) = \frac{(n_3 - n_1) \cdot (v_4 - v_2)}{2} - \bar{B}_{12} \quad (33)$$

$$c_{second,3}(f) = \frac{(n_4 - n_2) \cdot (v_3 - v_1)}{2} - \bar{B}_{21} \quad (34)$$

$$c_{second,4}(f) = \frac{(n_4 - n_2) \cdot (v_4 - v_2)}{2} - \bar{B}_{22} \quad (35)$$

The symbol $H_{ij}$ refers to the $ij$th component of the matrix $H$. Notice that the matrices $\bar{A}$ and $\bar{B}$ are not equal to the matrices $A$ and $B$ in Eq. (10). As the linear strain tensor $\varepsilon$ in Eq. (7) is relative to the identity matrix $I_2$, the matrices $A$ and $B$ in Eq. (10) can be regarded as the deformation trends of a unit square. Thus, a mapping needs to be established from the unit square to the original domain for $A$ and $B$, which is in the same way as setting up the mapping $\varphi_1: R^2 \to R^2$ in Section 3.1. By parameter transformation we have

$$\bar{A} = J^T \cdot A \cdot J \quad (36)$$

$$\bar{B} = J^T \cdot B \cdot J \quad (37)$$

The minimization of Eq. (28) does not prevent non-planarity, zigzagging or concave of quadrilaterals in the mesh. Thus, we added three regularizing terms $E_{coplanar}$, $E_{convex}$ and $E_{similar}$ as

$$E_{coplanar} = \omega_{coplanar} \sum_{f \in F} c_{coplanar}(f)^2 \qquad (38)$$

$$c_{convex}(f) = (v_3 - v_1) \cdot ((v_2 - v_1) \times (v_4 - v_1)) \qquad (39)$$

$$E_{convex} = \omega_{convex} \sum_{f \in F} c_{convex}(f)^2 \qquad (40)$$

$$c_{convex}(f) = \sum_{j=1}^{4} \theta_j - 2\pi \qquad (41)$$

$$E_{similar} = \omega_{similar} \sum_{f \in F} c_{similar}(f)^2 \qquad (42)$$

$$c_{similar}(f) = (v_3 + v_1 - v_2 - v_4)^2 - (v_3^M + v_1^M - v_2^M - v_4^M)^2 \qquad (43)$$

where $v_i$ is the $i$th vertex (in clockwise order) of the quadrilateral, $\theta_j$ is the $j$th inner angle of the quadrilateral, and $v_i^M$ is the $i$th vertex (in clockwise order) of the corresponding quadrilateral in initial mesh M. Then the final optimization goal can be modified as

$$E_{final} = E_{discrete} + E_{coplanar} + E_{convex} + E_{similar} \qquad (44)$$

We set $\omega_{coplanar} = \omega_{convex} = \omega_{similar} = 1$ to keep the coplanarity, convexity and the shapes of the quadrilaterals in the optimized mesh.

The solution of the non-linear least squares problem is realized with the LM algorithm [53], see Alg. 1. With the visualization of the solved results, the deformation simulation model of self-morphing bilayers is established.

**Alg. 1**

The Levenberg-Marquarelt method

| **Problem:** $\mathbf{x}^* = \text{argmin}_\mathbf{x}\{F(\mathbf{x})\}$ where $F(\mathbf{x}) = \frac{1}{2}\mathbf{f}(\mathbf{x})^T\mathbf{f}(\mathbf{x})$ |
|---|
| **begin** |
|   $k := 0$; $\mathbf{x} := \mathbf{x}_0$; $J(\mathbf{x}) := \partial \mathbf{f}(\mathbf{x})/\partial \mathbf{x}$; $\mathbf{A} := J(\mathbf{x})^T J(\mathbf{x})$ |
|   $\mathbf{g} := J(\mathbf{x})^T \mathbf{f}(\mathbf{x})$; $v := 2$; $\mu := \tau * \max\{A_{ij}\}$ |
|   $found := (\|\mathbf{g}\|_\infty \leq \varepsilon_1)$ |
|   **while** (**not** $found$) **and** ($k < k_{max}$) |
|     $k := k + 1$; Solve $(\mathbf{A} + \mu \mathbf{I})\mathbf{h}_{lm} = -\mathbf{g}$ |
|     **if** $\|\mathbf{h}_{lm}\| \leq \varepsilon_2(\|\mathbf{x}\| + \varepsilon_2)$ |

```
        found ≔ true
    else
        x ≔ x_new + h_lm
        ρ ≔ (F(x) − F(x_new))/(L(0) − L(h_lm))
        if ρ > 0
            x ≔ x_new;  A ≔ J(x)ᵀJ(x);  g ≔ J(x)ᵀf(x)
            found ≔ (‖g‖_∞ ≤ ε₁)
            μ ≔ μ ∗ max{1/3, 1 − (2ρ − 1)³};  ν ≔ 2
        else
            μ ≔ μ ∗ ν;  ν ≔ 2 ∗ ν
end
```

## 4 Performance of the deformation simulation model

To verify the performance of the proposed simulation model, we conducted some numerical experiments on bending rectangular bilayers. Bending is a common deformation pattern of bilayers that is caused by the one-directional strain difference of two layers. Two situations of bending were considered in our study: the small strain situation and the finite strain situation.

In the small strain situation, the strains of the two layers are assumed to be small quantities. The analysis of this situation, also known as pure bending, always follows the small strain assumptions: a) the second-order quantities of strain can be omitted, and b) there is no stretch of layers during deformation. The theoretical curvature of pure bending has been deduced by Timoshenko [54] and serves as a good reference for comparing the simulation results with the theoretical values. As such, pure-bending bilayers were an ideal case for verifying the performance of the proposed simulation model.

In the finite strain situation, the small strain assumptions no longer hold as the strains are not small enough. Timoshenko's theoretical formula for deformation curvature is also not applicable in this situation. Additionally, as we mentioned earlier, the simulation results of traditional FEM-based methods [40-43] are unreliable as they also rely on the small strain assumptions.

However, the proposed simulation model, which considers the extension of the layers (corresponding to the item $A$), offers an advantage in this situation. To illustrate this advantage, numerical experiments on bending bilayers in the finite strain situation were also conducted.

4.1 Numerical experiments on pure bending

For convenience, consider a simple bilayer rectangle with two layers having the same thickness ($t_1 = t_2 = t/2$) and Lamé constants ($\mu_1 = \mu_2 = \mu$). According to Timoshenko's work [54], the theoretical curvature of the pure-bending rectangle can be calculated as

$$\frac{1}{R_T} = \frac{3\Delta\varepsilon}{2t} \tag{45}$$

where $\Delta\varepsilon$ is the strain difference in the bilayer along the bending direction, $t$ is the thickness of the bilayer rectangle and $R_T$ is the curvature radius.

The energy items $A$ and $B$ in Eq. (10) are suitable for finite deformation situations. However, to comply with the small strain assumptions, these terms need to be modified to ensure that the simulation results are comparable to the theoretical values.

$$\omega_{first} = \frac{\mu t}{4}, \omega_{second} = \frac{\mu t^3}{4}, A = I_2, B = \frac{3\Delta\varepsilon}{2t} \tag{46}$$

The solver parameters in the numerical simulations of this study were listed as follows: $\tau = 0.001$, $\varepsilon_1 = 10\text{e-}8$, $\varepsilon_2 = 10\text{e-}8$, $k_{max} = 1000$. The bilayer rectangle had dimensions of 50mm×10mm×1mm and was meshed using 1mm×1mm quadrilateral elements. The strain difference $\Delta\varepsilon$ was varied from 0.01 to 0.05 with an interval of 0.01. All numerical experiments were conducted using both the proposed simulation model and the FEM-based model that we had previously proposed [42].

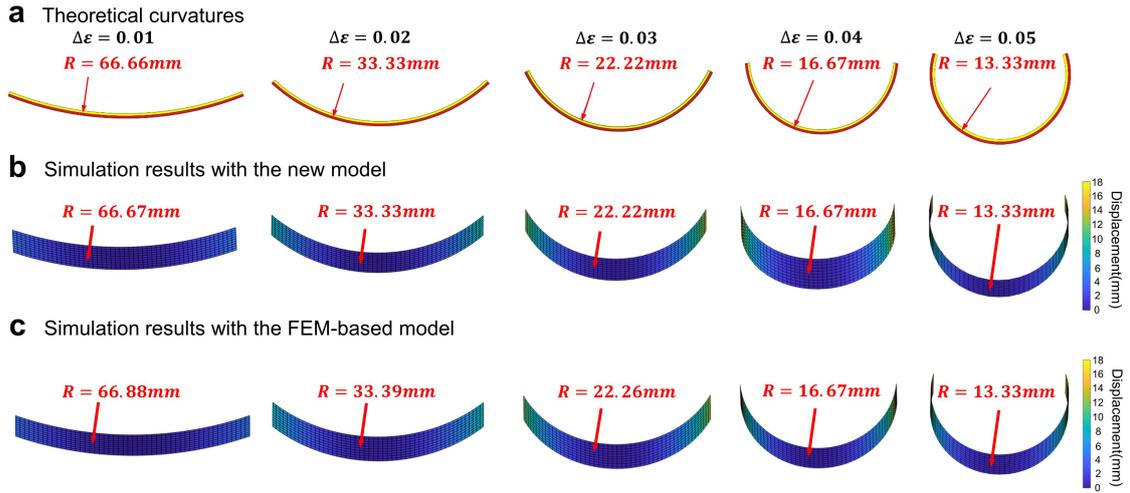

**Fig. 3.** (a) Theoretical curvatures of the bending rectangles with Δε from 0.01 to 0.05. (b) Simulation results of the bending rectangles with the checkerboard-based model in this work. (c) Simulation results of the bending rectangles with the FEM-based model in our previous work.

Table 1 and Figure 3 show the comparison between the theoretical values and the simulation results obtained from the proposed checkerboard-based simulation model and the FEM-based model. The simulation speeds of both models are also compared in Table 1. The results of the two models were comparable because they were both developed based on Matlab 2022b and were also run under the same environment. By comparison, the results of the proposed simulation model match the theoretical values better than the results of the FEM-based model, which demonstrates the accuracy of the new model. Moreover, the new simulation model is much faster compared to the FEM-based model due to the incorporation of the checkerboard-based DDG approach.

**Table 1**

Comparison of simulation accuracy and speeds on pure-bending deformation simulations between the proposed model and the FEM-based model

| Strain difference Δε | 0.01 | 0.02 | 0.03 | 0.04 | 0.05 |
| --- | --- | --- | --- | --- | --- |
| Theoretical curvature radii (mm) | 66.66 | 33.33 | 22.22 | 16.67 | 13.33 |
| Simulated curvature radii of the new model (mm) | 66.67 | 33.33 | 22.22 | 16.67 | 13.33 |
| Simulation time of new model (s) | 245 | 322 | 500 | 593 | 658 |
| Simulated curvature radii of the FEM-based model (mm) | 66.88 | 33.39 | 22.26 | 16.69 | 13.35 |
| Simulation time of the FEM-based model (s) | 14975 | 26134 | 40069 | 61886 | 90993 |

In addition, we conducted the numerical experiments again using the proposed simulation model, using two different element sizes (2mm×2mm and 0.5mm×0.5mm) to mesh the bilayer rectangle. The results (see Table 2) demonstrate that the proposed simulation model will converge to the theoretical values with refining meshes. Overall, these results demonstrate the effectiveness and efficiency of the new simulation model, allowing for powerful and efficient simulations in a range of scenarios.

**Table 2**

Performance of the proposed model on pure-bending deformation simulations under different mesh density

| Strain difference $\Delta\varepsilon$ | 0.01 | 0.02 | 0.03 | 0.04 | 0.05 |
|---|---|---|---|---|---|
| Theoretical curvature radii (mm) | 66.66 | 33.33 | 22.22 | 16.67 | 13.33 |
| Simulated curvature radii with 2mm×2mm elements (mm) | 66.67 | 33.34 | 22.24 | 16.68 | 13.35 |
| Simulated curvature radii with | 66.67 | 33.33 | 22.22 | 16.67 | 13.33 |
| 1mm×1mm elements (mm) | | | | | |
| Simulated curvature radii with 0.5mm×0.5mm elements (mm) | 66.67 | 33.33 | 22.22 | 16.67 | 13.33 |

4.2 Numerical experiments on bending bilayers with finite strains

Section 4.1 has already confirmed the accuracy of the proposed simulation model. However, to demonstrate the necessity and superiority of the model of considering the stretch of bilayers, we designed several additional numerical experiments with the simulation parameters shown in Table 5. In these experiments, the strain differences between the two layers were all 0.1. According to small strain theory, the curvature of these bending bilayers should be the same. We then compared the simulation results from our proposed model with the curvature radii calculated using Eq. (45) (see Table 3).

**Table 3**

Comparison between the actual curvature radii (simulated with our model) and the theoretical curvature radii (calculated using the small strain theory) of the bending bilayers with finite strains

| Top layer strain/bottom layer strain | 0.2/0.3 | 0.1/0.2 | 0/0.1 | -0.1/0 | -0.2/-0.1 | -0.3/-0.2 |
|---|---|---|---|---|---|---|
| Calculated curvature radii (mm) | 6.67 | 6.67 | 6.67 | 6.67 | 6.67 | 6.67 |
| Simulated curvature radii (mm) | 6.96 | 6.79 | 6.70 | 6.70 | 6.90 | 7.53 |

In addition, Fig. 4 highlights the significant differences between the simulation outcomes, further emphasizing the necessity of considering the stretch of bilayers in finite strain situations and the advantages of the proposed simulation model.

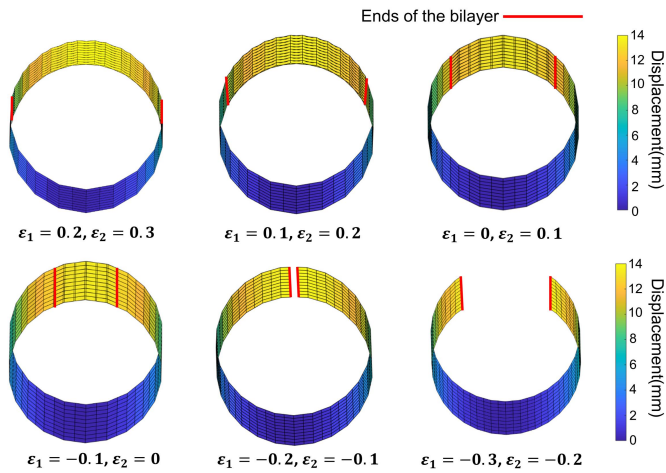

**Fig. 4.** Simulation results of the bending bilayers with finite strains

**5 Bio-inspired design and fabrication of 4D printed self-morphing bilayers**

The efficient simulation model exhibits its great application potential in the rapid verification of the design results of self-morphing bilayers. Inspired by some interesting creatures in nature, we designed some self-morphing bilayers, used the proposed simulation model to predict the growth morphology of the bilayers, and realized the fabrication of these surfaces by 4D printing technology. The simulation results were compared with the fabrication results, verifying the effect of the simulation model again.

The 4D printing technology used in this study is the FDM technology based on polylactic acid (PLA). The bilayer structures were printed by a FDM printer using PLA materials. Due to the thermosensitive characteristics of PLA, the printed bilayer will contain a prestress field at room temperature. When the printed structure is heated to $80°C$ (which is also the triggered condition for the deformation of the 4D printed PLA bilayers in our work), the prestress in the bilayer will be released and each layer will shrink along the printing path, which causes a mismatch between two layers and thus induce the deformation. The prestress of each layer in the bilayer is related to the printing speed of the layer. When the printing speed of the bottom layer is 300mm/min, the strain difference in the bilayer and the printing speed of the top layer has the relationship [55] as

$$\Delta\varepsilon = 3.593 \times 10^{-2} - 1.013 \times 10^{-4} v_t \tag{47}$$

where $v_t$ is the printing speed of the top layer. In this study, we chose the printing speed of the bottom layer in the bilayer as 300mm/min, then used this relationship to decide the printing speed of the top layer according to the strain difference we desired in the bilayer.

5.1 Four-petaled folding flowers

Four-petaled folding flowers are common in nature (like lilacs in Fig. 5a). Intuitively, we designed a bilayer structure, containing a square center and four petals around it. Fig. 5b shows the designed bilayer with a 16mm×16mm square center and a total size of 40mm×40mm×1mm. And the directions of strain differences were toward to the center due to the centripetal bending deformation of the petals, which also represent the printing directions of two layers. To show the different deformations under different strain differences, we designed the strain difference between two layers from 0.03 to 0.12 with an interval of 0.03. The printing speed of the top layer was set from 650mm/min to 1550mm/min with an interval of 300mm/min, while the printing speed of the bottom layer was 300mm/min, according to Eq. (47). The centripetal printing directions and the printing speed difference make the four petals of the printed bilayer bend toward center under heated, which deform like a lilac. Fig. 5c,d show the simulation result and the experiment result, coinciding well.

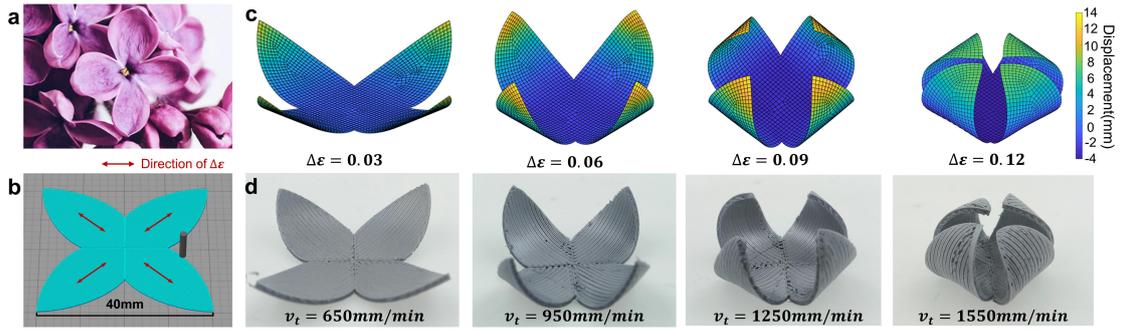

**Fig. 5.** (a) Lilac: a common four-petaled folding flower in nature. (b) The designed structure of the 4D printed self-morphing bilayer inspired by four-petaled folding flowers. (c) Simulation results of the flower-inspired self-morphing bilayers. (d) Experiment results of the flower-inspired self-morphing bilayers.

5.2 Spiral Grass

Some grasses in nature present interesting forms, such as twisting like springs, as shown in Fig. 6a. Here we imitated one kind of spiral grasses, designed the 4D printed bilayers as strip shapes, and changed the directions of strain difference to control the twisting degrees. Fig. 6b shows the designed shape of the bilayer (with a total size of 100mm×5mm×1mm). The angle $\gamma$ which indicated the printing direction of the bilayer (see Fig. 6b) was set from 30° to 75° with an interval of 15°. The strain difference between the two layers was set as 0.06 and thus the printing speed of the top layer was 950mm/min. The simulation and experiment results of the grass-inspired 4D printed PLA bilayers were shown in Fig. 6c,d.

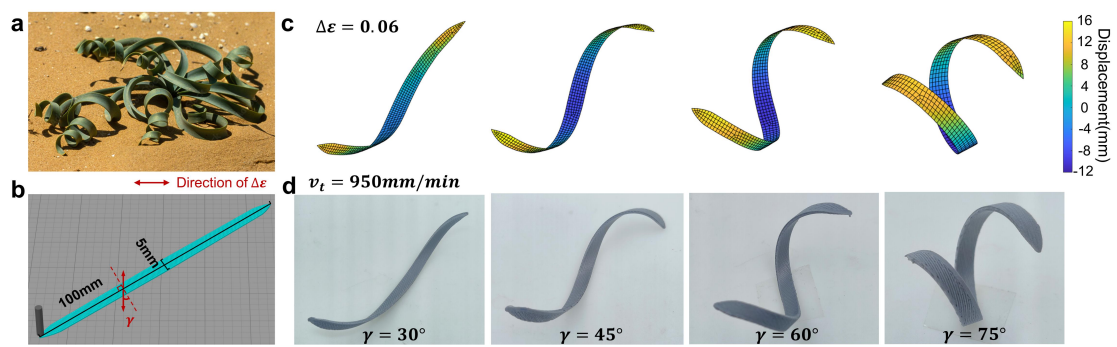

**Fig. 6.** (a) Albuca concordiana: a kind of spiral grasses in nature. (b) The designed structure of the 4D printed self-morphing bilayer inspired by spiral grasses. (c) Simulation results of the grass-inspired self-morphing bilayers. (d) Experiment results of the grass-inspired self-morphing bilayers.

5.3 Conical seashell

Inspired by the conical helicoid of conical seashells (see Fig. 7a), we try to explore how to generate a conical helicoid by strain mismatch. According to our previous work, a semi-annular shaped 4D printed PLA bilayer with concentric printing paths will deform to a cone-shape after triggered. Thus, we changed the circular printing paths in the semi-annular structure to the spiral printing paths using Archimedes spirals. Fig. 7b shows the designed shape of the obtained bilayer structure, where the point $O$ is the center of the Archimedes spirals. The parametric equation of Archimedes spiral can be written as

$$x = (\alpha + \beta\theta)\cos\theta \tag{48}$$

$$y = (\alpha + \beta\theta)\sin\theta \tag{49}$$

The outlines of the designed bilayer contain two Archimedes spirals, inner spiral ($\alpha = 20mm$ and $\beta = 1/18$) and outer spiral ($\alpha = 30mm$ and $\beta = 1/18$), with parameter $\theta$ starting from 0° to 180°. The printing directions (also the strain directions) of two layers are also Archimedes spirals with $20mm < \alpha < 30mm$, $\beta = 1/18$ and $\theta$ from 0° to 180° (see the red double arrow line in Fig. 7b). Fig. 7 c,d show the simulation and experiment results of the seashell-inspired 4D printed PLA bilayers, accompanied with the corresponding strain differences and the printing speeds.

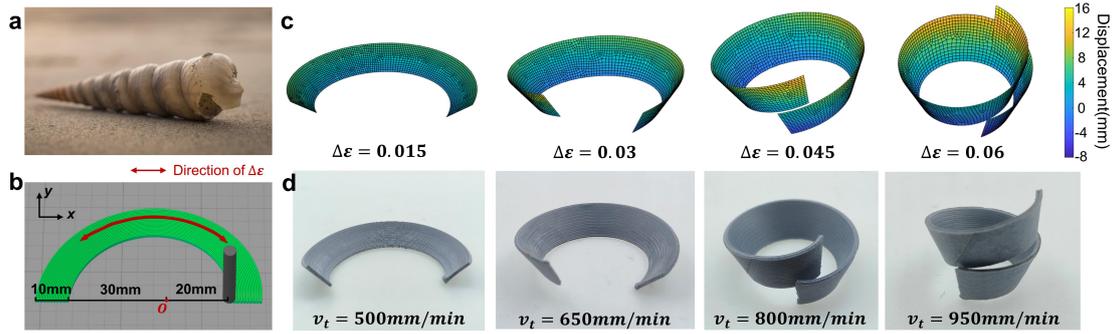

**Fig. 7.** (a) Turritella terebra: a kind of conical seashells in nature. (b) The designed structure of the 4D printed self-morphing bilayer inspired by conical seashells. (c) Simulation results of the seashell-inspired self-morphing bilayers. (d) Experiment results of the grass-inspired seashell-morphing bilayers.

In these design cases, the quadrilateral meshes in simulations were all unstructured. The simulation results still met the fabrication results well, showing that the simulation model has good accuracy for growth morphology prediction of irregular surfaces and performs well in self-morphing bilayer design.

## 6 Conclusion

In this work, we proposed an efficient and accurate deformation simulation model for self-morphing bilayers based on a newly proposed checkerboard-based DDG approach. The elastic energy function of a self-morphing bilayer was simplified and rearranged into a concise form, which intuitively illustrated the deformation trend of the bilayer during self-morphing. Moreover, we improved the checkerboard-based DDG approach to suit unstructured quadrilateral meshes, extending its use for more complex surfaces. We applied the improved DDG approach to discretize the deduced energy function, solve the minimization problem, and establish the deformation simulation model.

The performance of our proposed model is verified, including its accuracy and efficiency, by comparing it with a FEM-based simulation model from our previous work. Numerical results for large strain situations also demonstrated the superiority of the new model. We further applied the proposed model to design three bio-inspired bilayer structures, taking advantage of its high efficiency and reliability.

The introduced DDG approach based on checkerboard patterns is noteworthy for its concise forms and provable convergence. Our work is the first attempt to apply this approach in elastic energy minimization problems of plates. We believe our research will inspire other researchers to make good use of this approach, not just in computer graphics but also in mechanics problems.

One key advantage of DDG approaches relative to FEM approaches is the high computational efficiency. This feature also played an important role in our design process. Our future work will continue to use DDG approaches to capture the fundamental quantities of spatial surfaces, which will enable us to solve the inverse problem of obtaining the mismatch strain in a bilayer based on the desired deformation.

**Acknowledgment**

This work was supported by the National Natural Science Foundation of China (No. 52130501, No. 51935009, No. 52075479 and No. 52205288).